\magnification=\magstep1
\font\titre=cmbx10 scaled\magstep1
\font\soustitre=cmcsc10 scaled\magstep1

\font\aut=cmcsc10 scaled\magstep0
\font\abs=cmr8 scaled\magstep0
\def\BK{{\bf K}}
\def\GL{{\rm GL}}

\def\carre{\sqcap\mathrel{\mkern-22mu}\sqcup}
\def\BR{{\bf R}}
\centerline{\titre Commensurators of parabolic subgroups}
\centerline{\titre of Coxeter groups}
\bigskip
\centerline{\aut Luis Paris}
\bigskip
\centerline{Revised version}
\centerline{January 1996}
\bigskip\bigskip
\centerline{\vbox{\hsize=14truecm \abs
\noindent {\bf Abstract.}
Let $(W,S)$ be a Coxeter system, and let $X$ be a subset of
$S$. The subgroup of $W$ generated by $X$ is denoted by
$W_X$ and is called a parabolic subgroup.
We give the precise definition of the commensurator of a subgroup in a
group.
In particular, the commensurator of $W_X$ in $W$ is the subgroup
of $w$ in $W$ such that $wW_Xw^{-1}\cap W_X$ has finite index 
in both $W_X$ and $wW_Xw^{-1}$. The subgroup $W_X$
can be decomposed in the form $W_X = W_{X^0} \cdot
W_{X^\infty} \simeq W_{X^0} \times W_{X^\infty}$ where
$W_{X^0}$ is finite and all the irreducible components
of $W_{X^\infty}$ are infinite. Let $Y^\infty$ be the set
of $t$ in $S$ such that $m_{s,t}=2$ for all $s\in
X^\infty$. We prove that the commensurator of $W_X$ is
$W_{Y^\infty} \cdot W_{X^\infty} \simeq W_{Y^\infty}
\times W_{X^\infty}$. In particular, the commensurator of
a parabolic subgroup is a parabolic subgroup, and
$W_X$ is its own commensurator if and only if $X^0=
Y^\infty$.}}
\bigskip\bigskip
\vfootnote{}
{Mathematics Suject Classification: 20F55}
\centerline{\soustitre 1. Introduction}
\bigskip
Let $S$ be a finite set. A {\it Coxeter matrix} over $S$ is
a matrix $M= (m_{s,t})_{s,t\in S}$ indexed by the elements
of $S$ and satisfying\par
(a) $m_{s,s}=1$ if $s\in S$,\par
(b) $m_{s,t} = m_{t,s} \in \{ 2,3,4, \dots, +\infty \}$ if
$s,t\in S$ and $s\neq t$.\par\noindent
A Coxeter matrix $M=(m_{s,t})_{s,t\in S}$ is usually
represented by its {\it Coxeter graph} $\Gamma$. This is
defined by the following data.\par
(a) $S$ is the set of vertices of $\Gamma$.\par
(b) Two vertices $s,t\in S$ are joined by an edge if
$m_{s,t} \ge 3$.\par
(c) The edge joining two vertices $s,t\in S$ is labeled by
$m_{s,t}$ if $m_{s,t} \ge 4$.\par\noindent
The {\it Coxeter system} associated with $M$ (or with
$\Gamma$) is the pair $(W,S)$ where $W$ is the group
having the presentation
$$
W= \langle S\ |\ (st)^{m_{s,t}} =1\ {\rm if}\ m_{s,t} <
+\infty \rangle\ .
$$
The group $W$ is called the {\it Coxeter group} associated
with $M$. Given $X\subseteq S$, we write
\bigskip
\vbox{\halign{#\hfil\cr
$M_X=(m_{s,t})_{s,t\in X}$,\cr
$\Gamma_X$ the Coxeter graph which represents $M_X$,\cr
$W_X$ the subgroup of $W$ generated by $X$.\cr}}
\bigskip\noindent
The pair $(W_X,X)$ is the Coxeter system associated
with $M_X$ (see [Bo, Ch. IV, \S 1, $n^o$ 8]).
The group $W_X$ is called a {\it parabolic subgroup} of
the Coxeter system $(W,S)$. We assume that the reader is
familiar with the theory of Coxeter groups. We refer
to [Bo] and [Hu] for general expositions on the subject.
\par
For a group $G$ and for a subgroup $H$ of $G$, we denote by
$Z(G)$ the center of $G$, by $Z_G(H)$ the centralizer of
$H$ in $G$, by $N_G(H)$ the normalizer of $H$ in $G$, and
by $C_G(H)$ the commensurator of $H$ in $G$. Recall that
this is defined by
$$
C_G(H) = \{ g\in G\ ;\ H\cap (gHg^{-1})\ {\rm has\ finite\
index\ in\ both}\ H\ {\rm and}\ gHg^{-1} \}\ .
$$
Commensurators play an important role in representation
theory, especially in the study of induced representations.
For example, if a subgroup $H$ of $G$ is its own
commensurator, then any finite dimensional irreducible
representation of $H$ induces an irreducible
representation of $G$ (see [Ma]).
If $\BK$ is an infinite field and $P$ is a parabolic
subgroup of $\GL(n,\BK)$, then $P$ is its own
commensurator (see [BH]). A similar result is obviously
not true for Coxeter groups. Indeed, the commensurator of
a finite parabolic subgroup is the whole group $W$.
However, we prove in this paper that the commensurator of
a parabolic subgroup is always a parabolic subgroup
(Corollary 2.2),
and we give a criterion which decides whether a parabolic
subgroup is its own commensurator (Corollary 2.3).\par
The goal of this paper is to determine the commensurator
of a parabolic subgroup $W_X$ of a Coxeter system.
This subgroup can be decomposed in
the form $W_X = W_{X^0} \cdot W_{X^\infty} \simeq W_{X^0}
\times W_{X^\infty}$ where $W_{X^0}$ is finite and all the
irreducible components of $W_{X^\infty}$ are infinite. In
a first step (Proposition 2.4), we prove that the
commensurator of $W_X$ is the normalizer of
$W_{X^\infty}$. In a second step (Proposition 2.5), we
prove that the normalizer of $W_{X^\infty}$ is
$QZ_W(W_{X^\infty}) \cdot W_{X^\infty}$ where
$$
QZ_W(W_{X^\infty}) = \{ w\in W\ ;\ wX^\infty w^{-1} =
X^\infty \}
$$
is the {\it quasi-centralizer} of $W_{X^\infty}$. In a
third step (Proposition 2.6), we prove that\break
$QZ_W(W_{X^\infty})$ is $W_{Y^\infty}$ where $Y^\infty$ is
the set of $t$ in $S$ such that $m_{s,t}=2$ for all $s\in
X^\infty$. Finally, from Propositions 2.4, 2.5 and 2.6, we
deduce the following expression of the commensurator of
$W_X$ (Theorem 2.1).
$$
C_W(W_X) = W_{Y^\infty} \cdot W_{X^\infty} \simeq
W_{Y^\infty} \times W_{X^\infty} \ .
$$
\indent
We precisely state our results in Section 2, and we prove them in
Section 3.
\bigskip\bigskip
\centerline{\soustitre 2. Statements}
\bigskip
From now on, we fix a Coxeter system $(W,S)$.\par
Let $X$ be a subset of
$S$. Let $\Gamma_1, \dots, \Gamma_n$ be the connected
components of $\Gamma_X$ and, for $i\in \{ 1,\dots, n\}$,
let $X_i$ be the set of vertices of $\Gamma_i$. The group
$W_{X_i}$ is called an {\it irreducible component} of
$W_X$. It is clear that
$$
W_X = W_{X_1} \cdot \dots \cdot W_{X_n} \simeq
W_{X_1} \times \dots \times W_{X_n}\ .
$$
We assume that $W_{X_i}$ is finite if $i=1, \dots, r$, and
that $W_{X_i}$ is infinite if $i=r+1, \dots, n$. We set
$$\eqalign{
&X^0= X_1 \cup \dots \cup X_r\ ,\cr
&X^\infty = X_{r+1} \cup \dots \cup X_n\ .\cr}
$$
Then
$$
W_X = W_{X^0} \cdot W_{X^\infty}\
\simeq W_{X^0} \times W_{X^\infty}\ ,
$$
the group $W_{X^0}$ is finite, and all the irreducible
components of $W_{X^\infty}$ are infinite.
\bigskip\noindent
{\aut Theorem 2.1.} {\it Let $X$ be a subset of $S$. Then
$$
C_W(W_X) = W_{Y^\infty} \cdot W_{X^\infty} =
W_{Y^\infty \cup X^\infty} \simeq
W_{Y^\infty} \times W_{X^\infty}
$$
where
$$
Y^\infty = \{ t\in S\ ;\ m_{s,t} =2\ {\it for\ all}\
s\in X^\infty \}\ .
$$}
\bigskip\noindent
{\aut Corollary 2.2.} {\it The commensurator of a
parabolic subgroup of $(W,S)$ is a parabolic subgroup.}
\bigskip\noindent
{\aut Corollary 2.3.} {\it Let $X$ be a subset of $S$.
Then $W_X$ is its own commensurator if and only if
$X^0$ is the set of $t\in S$ such that $m_{s,t}=2$ for
all $s\in X^\infty$.}
\bigskip
Theorem 2.1 is a direct consequence of the following
Propositions 2.4, 2.5, and 2.6.
\bigskip\noindent
{\aut Proposition 2.4.} {\it Let $X$ be a subset of $S$. Then
the commensurator of
$W_X$ in $W$ is equal to the normalizer of $W_{X^\infty}$
in $W$.}
\bigskip
We define the {\it quasi-center} of $(W,S)$
to be
$$
QZ(W,S) = \{ w\in W\ ;\ wSw^{-1} = S\}\ .
$$
Similarly,
we define the {\it quasi-centralizer} of a parabolic
subgroup $W_X$ of $(W,S)$ to be
$$
QZ_W(W_X) = \{ w\in W\ ;\ wXw^{-1} =X \}\ .
$$
\bigskip\noindent
{\aut Proposition 2.5.} {\it Let $X$ be a subset of $S$. Then
$$
N_W(W_X) = QZ_W(W_X) \cdot W_X\ .
$$
Moreover,
$$
QZ_W(W_X) \cap W_X = QZ(W_X,X)\ .
$$}
\bigskip\noindent
{\aut Proposition 2.6.} {\it Let $X$ be a subset of $S$
such that all the
irreducible components of $W_X$ are infinite (i.e.
$X=X^\infty$). Let
$$
Y=\{ t\in S\ ;\ m_{s,t} =2\ {\it for\ all}\ s\in X \}\ .
$$
Then the quasi-centralizer of $W_X$ is equal to
$W_Y$.}
\bigskip
Proposition 2.4 is a consequence of [So, Lemma 2].
Proposition 2.5 is stated in [Ho] for finite type Coxeter systems
(see also [Kr, Ch. 3]).
Moreover, its proof is quite simple. Proposition 2.6 is a
consequence of [De, Prop. 5.5].
\bigskip\bigskip
\centerline{\soustitre 3. Proofs}
\bigskip
First, we state in Lemmas 3.1 and 3.2 some well-known
facts that will be required later.
Recall that each $w\in W$ can be written $w=s_1 \dots s_r$
where $s_i\in S$ for all $i\in \{ 1, \dots, r\}$. If $r$
is as small as possible, then $r$ is called the {\it
length} of $w$ and is denoted by $l(w)$.
\bigskip\noindent
{\aut Lemma 3.1} ({\aut Bourbaki} [Bo, Ch. IV, \S 1, Ex. 3]).
{\it Let $X$ and $X'$ be two subsets of $S$.\par
(i) Let $w\in W$. There is a unique element $v$ of minimal
length in $W_XwW_{X'}$. Moreover, each $w'\in W_X w
W_{X'}$ can be written as $w'= uvu'$ where $u\in W_X$, $u'\in
W_{X'}$, and $l(w') = l(u) + l(v) + l(u')$.\par
An element $v$ is called $(X,X')$-reduced if it is of
minimal length in $W_X v W_{X'}$.\par
(ii) If an element $v$ is $(X, \emptyset)$-reduced, then
$l(uv) = l(u) + l(v)$ for all $u\in W_X$.\par
(iii) If an element $v$ is $(\emptyset,X')$-reduced, then
$l(vu') = l(v) + l(u')$ for all $u'\in W_{X'}$.\par
(iv) An element $v$ is $(X, \emptyset)$-reduced if and only
if $l(sv)>l(v)$ for all $s\in X$.\par
(v) An element $v$ is $(\emptyset, X')$-reduced if and only
if $l(vs') > l(v)$ for all $s'\in X'$.\par
(vi) An element $v$ is $(X,X')$-reduced if and only if it is
both $(X,\emptyset)$-reduced and $(\emptyset,
X')$-reduced.}
\bigskip\noindent
{\aut Lemma 3.2} ({\aut Bourbaki} [Bo, Ch. IV, \S 1, Ex. 22]).
{\it Let $w_0$ be an element of $W$. The following
statements are equivalent.\par
(1) $l(sw_0)<l(w_0)$ for all $s\in S$.\par
(2) $l(w_0s)<l(w_0)$ for all $s\in S$.\par
(3) $l(ww_0) = l(w_0) - l(w)$ for all $w\in W$.\par
(4) $l(w_0w) = l(w_0) - l(w)$ for all $w\in W$.\par\noindent
Such an element is unique and exists if and only if $W$ is
finite. Then it is the unique element of maximal
length in W. Moreover, $w_0^2=1$ and $w_0 S w_0 = S$.}
\bigskip
The following proposition is the key of the proof of
Proposition 2.4.
\bigskip\noindent
{\aut Proposition 3.3} ({\aut Solomon} [So, Lemma 2]).
{\it Let $X$ and $X'$ be two subsets of $S$, and let $v$ be a
$(X,X')$-reduced element of $W$. Then
$$
W_X\cap (vW_{X'}v^{-1})=W_Y
$$
where $Y=(vX'v^{-1})\cap X$.}
\bigskip\noindent
{\aut Corollary 3.4.}
{\it Let $X$ and $X'$ be two subsets of $S$, and let
$w$ be an element of $W$. We write $w=u_0vu_0'$ where 
$u_0\in W_X$, $u_0'\in W_{X'}$, and $v$ is $(X,X')$-reduced.
Then
$$
W_X\cap (wW_{X'}w^{-1}) = u_0 W_Y u_0^{-1}
$$
where $Y=(vX'v^{-1})\cap X$.}
\bigskip\noindent
{\it Proof.}
$$\eqalign{
W_X \cap (w W_{X'} w^{-1})
=\ &W_X \cap (u_0vu_0' W_{X'} u_0'^{-1} v^{-1} u_0^{-1})\cr
=\ &W_X \cap (u_0v W_{X'} v^{-1} u_0^{-1})\cr
=\ &u_0 ((u_0^{-1} W_X u_0) \cap (v W_{X'} v^{-1}))
u_0^{-1}\cr
=\ &u_0 (W_X \cap (v W_{X'} v^{-1})) u_0^{-1}\cr
=\ &u_0 W_Y u_0^{-1}\ .\quad\carre\cr}
$$
\bigskip\noindent
{\it Proof of Proposition 2.4.} Let $w\in N_W(W_{X^\infty})$.
Then
$$
W_{X^\infty} = w W_{X^\infty} w^{-1} \subseteq W_X \cap (w
W_X w^{-1})\ ,
$$
the group $W_{X^\infty}$ has finite index in $W_X$, and the
group $w W_{X^\infty} w^{-1}$ has finite index in $w W_X
w^{-1}$. Thus $W_X \cap (w W_X w^{-1})$ has finite index
in both $W_X$ and $w W_X w^{-1}$. This shows that $N_W
(W_{X^\infty}) \subseteq C_W(W_X)$.\par
Let $w\in C_W(W_X)$. We write $w=u_0vu_0'$ where $u_0,u_0'
\in W_X$ and $v$ is $(X,X)$-reduced.
By Corollary 3.4,
$$
W_X \cap (w W_X w^{-1}) = u_0 W_Y u_0^{-1}
$$
where $Y=(vXv^{-1})\cap X$. Let $Y^0 = Y\cap X^0$, and let
$Y^\infty = Y \cap X^\infty$. For a group $G$ and for a
subgroup $H$ of $G$, we denote by $|G:H|$ the index of $H$
in $G$. Then
$$\displaylines{
|W_X : W_X \cap (w W_X w^{-1})| = |W_X: u_0 W_Y u_0^{-1}| =
|W_X: W_Y|\cr
= |W_{X^0}: W_{Y^0}| \cdot |W_{X^\infty}:
W_{Y^\infty}|\ .\cr}
$$
If $Y^\infty \neq X^\infty$, then, by [De, Prop. 4.2],
$W_{Y^\infty}$ has infinite index in $W_{X^\infty}$, thus
$W_X \cap (w W_X w^{-1})$ has infinite index in $W_X$, too.
This is not the case, thus $Y^\infty = X^\infty$. Let
$\Gamma_1, \dots, \Gamma_n$ be the connected components of
$\Gamma_X$, and, for $i=1, \dots, n$, let $X_i$ be the set
of vertices of $\Gamma_i$. We assume that $X^0 = X_1 \cup
\dots \cup X_r$ and that $X^\infty = X_{r+1} \cup \dots
\cup X_n$. Let $i\in \{ r+1, \dots, n\}$. Then
$$
v^{-1} X_i v \subseteq v^{-1} X^\infty v = v^{-1} Y^\infty
v \subseteq v^{-1} Y v \subseteq X\ .
$$
Thus there exists $j\in \{ 1, \dots, r, r+1, \dots, n\}$
such that $v^{-1} X_i v \subseteq X_j$. The group
$W_{X_i}$ is infinite and $v^{-1} W_{X_i} v \subseteq
W_{X_j}$, thus $W_{X_j}$ is infinite, therefore $j\in \{
r+1, \dots, n\}$. This shows that $v^{-1} X^\infty v
\subseteq X^\infty$, thus $v X^\infty v^{-1} = X^\infty$,
therefore $v W_{X^\infty} v^{-1} = W_{X^\infty}$. On the other hand,
since $W_X =W_{X^0} \cdot W_{X^\infty} \simeq
W_{X^0} \times W_{X^\infty}$, we have
$u W_{X^\infty} u^{-1} = W_{X^\infty}$ for all $u\in W_X$.
So,
$$
w W_{X^\infty} w^{-1} = u_0 v u_0' W_{X^\infty} u_0'^{-1}
v^{-1} u_0^{-1} = u_0 v W_{X^\infty} v^{-1} u_0^{-1} = u_0
W_{X^\infty} u_0^{-1} = W_{X^\infty}\ .
$$
This shows that $C_W(W_X) \subseteq N_W (W_{X^\infty})$.
$\carre$
\bigskip\noindent
{\it Proof of Proposition 2.5.} The inclusion
$$
QZ_W (W_X) \cdot W_X \subseteq N_W(W_X)
$$
is obvious.\par
Let $w\in N_W(W_X)$. We write $w=vu$ where $u\in W_X$, $v$
is $(\emptyset, X)$-reduced, and $l(w) = l(v) + l(u)$. We
have
$$
w W_X w^{-1} = v W_X v^{-1} = W_X\ .
$$
The element $v$ is of minimal length in $vW_X = W_X v$,
thus $v$ is also $(X, \emptyset)$-reduced. If $s\in X$,
then, by Lemma 3.1,
$$\eqalign{
&l(v)+1 = l(vs) = l(vsv^{-1}v) = l(vsv^{-1}) + l(v)\cr
\Rightarrow\ &l(vsv^{-1}) =1\cr
\Rightarrow\ &vsv^{-1} \in W_X \cap S = X\ .\cr}
$$
So, $vXv^{-1} \subseteq X$, thus $vXv^{-1} =X$, therefore
$v\in QZ_W(W_X)$. This shows that $N_W(W_X) \subseteq
QZ_W(W_X) \cdot W_X$.\par
The equality
$$
QZ_W(W_X) \cap W_X = QZ(W_X, X)
$$
is obvious. $\carre$
\bigskip
Before proving Proposition 2.6, we recall some facts on root
systems. Let $V$ be a real vector space having a basis $\{
e_s; s\in S \}$ in one-to-one correspondence with $S$. Let
$B$ be the symmetric bilinear form on $V$ defined by
$$
B(e_s,e_t) = \left\{
\eqalign{
&- \cos (\pi /m_{s,t}) \quad {\rm if}\ m_{s,t} <+\infty\cr
&-1\quad {\rm if}\ m_{s,t}=+\infty\cr}
\right.
$$
There is an action of $W$ on $V$ defined by
$$
s(x) = x- 2B(x,e_s) e_s
$$
if $s\in S$ and $x\in V$. This action is called the {\it
canonical representation} of $(W,S)$. The {\it root
system} $\Phi$ of $(W,S)$ is the collection of all vectors
$w(e_s)$ where $w\in W$ and $s\in S$. By [Bo, Ch. V, \S 4,
Ex. 8], every root $\alpha$ can be uniquely written in the
form
$$
\alpha = \sum_{s\in S} a_se_s\quad (a_s \in \BR)
$$
where either all $a_s$ are positive, or all $a_s$ are
negative. We call $\alpha$ {\it positive} and write
$\alpha >0$ if $a_s\ge 0$ for all $s\in S$. We call
$\alpha$ {\it negative} and write $\alpha <0$ if $a_s \le
0$ for all $s\in S$.
\bigskip\noindent
{\aut Proposition 3.5} ({\aut Deodhar} [De, Prop. 3.1]).
{\it Let
$$
T= \{ wsw^{-1}\ ;\ w\in W\ {\it and}\ s\in S \}\ ,
$$
and let $\Phi^+$ be the set of positive roots. For $\alpha
= w(e_s)$, we write $r_\alpha = wsw^{-1}$. Then the
function $\Phi^+ \to T$ $(\alpha \mapsto r_\alpha)$ is
well-defined and bijective.}
\bigskip\noindent
{\aut Proposition 3.6} ({\aut Deodhar} [De, Prop. 2.2]).
{\it Let $w\in W$, and let $s\in S$. Then $l(ws) > l(w)$
if and only if $w(e_s) >0$.}
\bigskip
For a subset $X$ of $S$, we write
$$
E_X = \{ e_s\ ;\ s\in X \}\ .
$$
The following lemma is an easy consequence of
Propositions 3.5 and 3.6.
\bigskip\noindent
{\aut Lemma 3.7.} {\it Let $X$ and $X'$ be two subsets
of $S$, and let $w$ be an element of $W$. The following
statements are equivalent.\par
(1) $w(E_X) = E_{X'}$.\par
(2) $wXw^{-1} = X'$ and $l(ws) > l(w)$ for all $s\in X$.}
\bigskip
For $X\subseteq S$ such that $W_X$ is finite, we denote by
$w_X$ the unique element of maximal length in $W_X$.\par
Let $X$ be a subset of $S$, and let $t$ be an element of $S
\setminus X$. Let $\Gamma_0$ be the connected component of
$\Gamma_{ \{t\} \cup X}$ containing $t$, and let $Y_0$ be
the set of vertices of $\Gamma_0$. We say that $t$ is {\it
$X$-admissible} if $W_{Y_0}$ is finite. In that case, we
write
$$
c(t,X) = w_{Y_0} w_{X_0}
$$
where $X_0 = Y_0 \setminus \{t\}$.
It is the element of minimal length in $w_{Y_0}W_X$. In
particular, $c(t,X)$ is $(\emptyset, X)$-reduced.
By Lemma 3.2 and
Lemma 3.7, there exists a subset $X'$ of $\{t\} \cup
X$ such that
$$
c(t,X) (E_X) = E_{X'}\ .
$$
If $X=X^\infty$, then $t$ is $X$-admissible if and only if
$m_{s,t}=2$ for all $s\in X$. In that case, $c(t,X)=t$ and
$c(t,X) (E_X) = E_X$.
\bigskip\noindent
{\aut Proposition 3.8} ({\aut Deodhar} [De, Prop. 5.5]). {\it
Let $X$ and $X'$ be two subsets of $S$, and let $w$ be an
element of $W$. If $w( E_X) = E_{X'}$, then there
exist sequences\par
$X_0=X, X_1, \dots, X_n=X'$ of subsets of $S$,\par
$t_0, t_1, \dots, t_{n-1}$ of elements of $S$,\par\noindent
such that\par
(a) $t_i\in S\setminus X_i$ and $t_i$ is $X_i$-admissible
($i=0,1, \dots, n-1$),\par
(b) $c(t_i,X_i) (E_{X_i}) = E_{X_{i+1}}$ ($i=0,1,
\dots, n-1$),\par
(c) $w= c(t_{n-1}, X_{n-1}) \dots c(t_1, X_1) c(t_0,
X_0)$.}
\bigskip
The following Lemmas 3.9 and 3.10 are preliminary results
 to the proof of Proposition 2.6.
\bigskip\noindent
{\aut Lemma 3.9.} {\it Let $X$ and $X'$ be two subsets of
 $S$, and let $w$ be an element of $W$. If $w X w^{-1} =
 X'$, then $w$ can be written $w =vu$ where $u\in
 QZ(W_X,X)$, $vXv^{-1} = X'$, and $l(vs) > l(v)$ for all
 $s\in X$.}
\bigskip\noindent
{\it Proof.} We write $w= vu$ where $u\in W_X$, $v$ is
 $(\emptyset, X)$-reduced, and $l(w) = l(v) + l(u)$. We
 have
$$
w W_X w^{-1} = v W_X v^{-1} = W_{X'}\ .
$$
The element $v$ is of minimal length in $vW_X = W_{X'} v$,
 thus $v$ is also $(X', \emptyset)$-reduced. If $s\in X$,
 then, by Lemma 3.1,
$$\eqalign{
&l(v)+1 = l(vs) = l(vsv^{-1}v) = l(vsv^{-1}) + l(v)\cr
\Rightarrow\ &l(vsv^{-1}) =1\cr
\Rightarrow\ &vsv^{-1} \in W_{X'} \cap S =X'\ .\cr}
$$
So, $vXv^{-1} \subseteq X'$. Similarly, $v^{-1} X' v
 \subseteq X$. Thus $vX v^{-1} = X'$.\par
Since $v$ is $(\emptyset , X)$-reduced, by Lemma 3.1,
 $l(vs) > l(v)$ for all $s\in X$.\par
Finally,
$$\eqalign{
&wXw^{-1} = vu X u^{-1}v^{-1} = X'\cr
\Rightarrow\ &uXu^{-1} = v^{-1} X'v = X\cr}
$$
thus $u\in QZ(W_X,X)$. $\carre$
\bigskip\noindent
{\aut Lemma 3.10} ({\aut Bourbaki} [Bo, Ch. V, \S 4, Ex. 3]).
{\it We suppose that $(W,S)$ is irreducible.\par
(i) If $W$ is finite, then $QZ(W,S) = \{1, w_0\}$, where
$w_0$ is the unique element of maximal length in $W$.\par
(ii) If $W$ is infinite, then $QZ(W,S) = \{ 1\}$.}
\bigskip\noindent
{\it Proof of Proposition 2.6.} The inclusion
$$
W_Y \subseteq QZ_W(W_X)
$$
is obvious.\par
Let $w\in QZ_W(W_X)$. By Lemma 3.9, $w$ can be written
$w=vu$ where $u\in$
\break
$QZ(W_X,X)$, $vXv^{-1} = X$, and
$l(vs)>l(v)$ for all $s\in X$. Since $X=X^\infty$, by
Lemma 3.10, $QZ(W_X,X) = \{1\}$, thus $u=1$. By Lemma
3.7, $v(E_X) = E_X$. By Proposition 3.8, there exist
sequences\par
$X=X_0,X_1, \dots, X_n=X$ of subsets of $S$,\par
$t_0, t_1, \dots, t_{n-1}$ of elements of $S$,\par\noindent
such that\par
(a) $t_i\in S \setminus X_i$ and $t_i$ is $X_i$-admissible
($i=0,1, \dots, n-1$),\par
(b) $c(t_i,X_i) (E_{X_i}) = E_{X_{i+1}}$ ($i=0,1,
\dots, n-1$),\par
(c) $v= c(t_{n-1},X_{n-1}) \dots c(t_1,X_1)
c(t_0,X_0)$.\par\noindent
Since $X=X^\infty$, if $X_i=X$, then $m_{t_i,s}=2$ for all
$s\in X$ (namely, $t_i\in Y$), $c(t_i,X_i) = t_i$, and
$X_{i+1} = X$ (since $t_i (E_X) = E_X$). Since
$X_0=X$, it follows that $c(t_i,X_i) = t_i \in Y$ for all
$i=0,1,\dots, n-1$. Thus
$$
w=v= t_{n-1} \dots t_1t_0 \in W_Y\ .
$$
This shows that $QZ_W(W_X) \subseteq W_Y$. $\carre$
\bigskip\bigskip
\centerline{\soustitre References}
\bigskip
\hbox{
\vtop{\hsize=1.3 truecm \noindent
[Bo]
\hfill}
\vtop{\hsize=14.7 truecm \noindent
{\aut N. Bourbaki},
``Groupes et alg\`ebres de Lie, Chapitres IV--VI'',
Hermann, Paris, 1968.
}}
\medskip
\hbox{
\vtop{\hsize=1.3 truecm \noindent
[Br]
\hfill}
\vtop{\hsize=14.7 truecm \noindent
{\aut K.S. Brown},
``Buildings'',
Springer-Verlag, New York, 1989.
}}
\medskip
\hbox{
\vtop{\hsize=1.3 truecm \noindent
[BH]
\hfill}
\vtop{\hsize=14.7 truecm \noindent
{\aut M. Burger} and {\aut P. de la Harpe},
{\it Irreducible representations of discrete groups},
in preparation.
}}
\medskip
\hbox{
\vtop{\hsize=1.3 truecm \noindent
[De]
\hfill}
\vtop{\hsize=14.7 truecm \noindent
{\aut V.V. Deodhar},
{\it On the root system of a Coxeter group},
Comm. Algebra {\bf 10} (1982), 611-630.
}}
\medskip
\hbox{
\vtop{\hsize=1.3 truecm \noindent
[Ho]
\hfill}
\vtop{\hsize=14.7 truecm \noindent
{\aut R.B. Howlett},
{\it Normalizers of parabolic subgroups of reflection groups},
J. London Math. Soc. (2) {\bf 21} (1980), 62--80.
}}
\medskip
\hbox{
\vtop{\hsize=1.3 truecm \noindent
[Hu]
\hfill}
\vtop{\hsize=14.7 truecm \noindent
{\aut J.E. Humphreys},
``Reflection groups and Coxeter groups'',
Cambridge studies in advanced mathematics, vol. 29,
Cambridge University Press, 1990.
}}
\medskip
\hbox{
\vtop{\hsize=1.3 truecm \noindent
[Kr]
\hfill}
\vtop{\hsize=14.7 truecm \noindent
{\aut D. Krammer},
``The conjugacy problem for Coxeter groups'',
Ph. D. Thesis, Utrecht, 1994.
}}
\medskip
\hbox{
\vtop{\hsize=1.3 truecm \noindent
[Ma]
\hfill}
\vtop{\hsize=14.7 truecm \noindent
{\aut G.W. Mackey},
``The theory of unitary group representations'',
The University of Chicago Press, 1976.
}}
\medskip
\hbox{
\vtop{\hsize=1.3 truecm \noindent
[So]
\hfill}
\vtop{\hsize=14.7 truecm \noindent
{\aut L. Solomon},
{\it A Mackey formula in the group ring of a Coxeter group},
J. Algebra {\bf 41} (1976), 255--268.
}}
\bigskip\bigskip
\vbox{\halign{#\hfil\cr
{\it Author's address:}\cr
\noalign{\smallskip}
Universit\'e de Bourgogne\cr
U.R.A. 755 -- Laboratoire de Topologie\cr
D\'epartement de Math\'ematiques\cr
B.P. 138\cr
21004 Dijon Cedex\cr
FRANCE\cr}}
\bigskip
\vbox{\halign{#\hfil\cr
{\it E-mail address:}\cr
\noalign{\smallskip}
lparis@satie.u-bourgogne.fr\cr}}
\end